\documentclass[12pt]{amsart}

\usepackage{amsfonts}
\usepackage{amsmath}
\usepackage{amsthm}

\newtheorem{theorem}{Theorem}[section]
\newtheorem{lemma}[theorem]{Lemma}
\newtheorem{coro}[theorem]{Corollary}

\def\oti{\tilde{\omega}^{-1}} 
\def\ot{\tilde{\omega}}
\def\etheta{\epsilon(\tilde{\omega}^{-1}\frac{\overline{\theta}}{\theta},\psi_0)}
\def\rtheta{r_{\theta}}
\def\nsi0{n(\psi _0)}
\def\otitheta{\tilde{\omega}^{-1}\frac{\overline{\theta}}{\theta}}
\def\tsbt{\frac{\overline{\theta}}{\theta}}
\def\tbts{\frac{\theta}{\overline{\theta}}}
\def\eoti{\epsilon(\tilde{\omega}^{-1},\psi_0)}
\def\eli{\epsilon(\lambda^{-1},\psi_0)}
\def\eltheta{\epsilon(\lambda^{-1}\frac{\overline{\theta}}{\theta},\psi_0)}
\def\elitheta{\epsilon(\lambda'^{-1}\frac{\overline{\theta}}{\theta},\psi_0)}

\def\gl2f{GL(2,F)}
\def\etsbt{\epsilon(\frac{\overline{\theta}}{\theta},\psi_0)}
\def \rem {\noindent{\bf Remark: }}

\begin{document}

\title[On Counting Twists of a Character ...]{On Counting Twists of a Character Appearing in its Associated Weil Representation}

\author{K. Vishnu Namboothiri}
\address{Department of Mathematics and Statistics, University of Hyder-
abad, Hyderabad - 500 046, India \& Department of Collegiate Education, Government of Kerala, Kerala, India}
\email{kvnamboothiri@gmail.com}

\begin{abstract}
Consider an irreducible, admissible representation $\pi$ of
GL(2,$F$) whose restriction to GL(2,$F)^+$ breaks up as a sum of two
irreducible representations $\pi_+ + \pi_-$. If $\pi=r_{\theta}$,
the Weil representation of GL(2,$F$) attached to a character
$\theta$ of $K^*$ which does not factor through the norm map from
$K$ to $F$, then $\chi\in \widehat{K^*}$ with $(\chi . \theta
^{-1})\vert _{ F^{ * }}=\omega _{ {K/F}}$ occurs in ${r_{\theta}}_+$
if and only if $\epsilon(\theta\chi^{-1},\psi_0)=\epsilon(\overline
\theta\chi^{-1},\psi_0)=1$ and in ${r_{\theta}}_-$ if and only if
both the epsilon factors are $-1$. But given a conductor $n$, can we
say precisely how many such $\chi$ will appear in $\pi$? We
calculate the number of such characters at each given conductor $n$
in this work.
\end{abstract}

\maketitle

\section{Introduction}

 Let $F$ be a nonarchimedean local field of
characteristic not two and $K$ a separable quadratic extension. Then
if $K=F(x_0)$ with $x_0$ an element of $K^*$ whose trace to $F$ is 0
we have an embedding of $K^*$ into GL(2,$F$) given by \[a+bx_0
\mapsto \left[\begin{array}{cc}a & bx^2_0\\b &
a\end{array}\right].\] .

Let GL(2,$F)^+$ be the subgroup of index 2 in GL(2,$F$) consisting
of those matrices whose determinant is in $N_{K/F}(K^*)$ where
$N_{K/F}$ is the usual norm map from $K$ to $F$. In \cite{prasad} D.
Prasad  considered  irreducible, admissible  representations $\pi$
of GL(2,$F$) whose restriction to GL(2,$F)^+$ breaks up as a sum of
two irreducible representations $\pi_+ + \pi_-$ . There he gave a
characterization of characters $\chi$ of $K^*$ occurring  in the
restriction of $\pi$ to $K^*$. It is immediate that if a character
$\chi$ occurs in such a restriction then $\chi|_{F^*}$ must be the
central character of $\pi$.  Hence if $\pi$ is supercuspidal then
$\pi=r_{\theta}$, the Weil representation of GL(2,$F$) attached to a
character $\theta$ of $K^*$ which does not factor through the norm
map from $K$ to $F$. Prasad showed that $\chi$ occurs in
${r_{\theta}}_+$ if and only if
$\epsilon(\theta\chi^{-1},\psi_0)=\epsilon(\overline
\theta\chi^{-1},\psi_0)=1$ and in ${r_{\theta}}_-$ if and only if
both the epsilon factors are $-1$ ($\overline \theta$ is the Galois
conjugate of $\theta$). What he proved exactly is the following:

\begin{theorem}\label{prasad}
Let $\rtheta$ be an irreducible admissible representation of
$GL(2,F)$ associated to a regular character $\theta$ of $K^{ *}$.
Fix embeddings of $K^{
*}$ in $GL(2,F)^{ +}$ and in ${D^{ {*+}}}_{ F}$, and
choose a nontrivial additive character $\psi$ of $F$, and an element
$x_{ 0}$ of $K^{ *}$ with $tr(x_{ 0})=0$. Then the representation
${\rtheta}$ of $GL(2,F)$ decomposes as ${\rtheta} = {\rtheta} _+
\oplus {\rtheta} _- $ when restricted to $GL(2,F)^+$ and the
representation ${\rtheta} '$ of ${D^{ *}}_{ F}$ decomposes as
${\rtheta} ' = {\rtheta} '_+ \oplus {\rtheta} '_- $ when restricted
to ${D^{ {*+}}}_{ F}$, such that for a character $\chi$ of $K^{
*}$ with $(\chi . \theta ^{-1})\vert _{ F^{ * }}=\omega
_{ {K/F}}$, $\chi$ appears in ${\rtheta} _+$ if and only if
$\epsilon (\theta \chi ^{-1},\psi_{ 0})=\epsilon (\overline{\theta}
\chi ^{-1},\psi_{ 0})=1$, $\chi$ appears in ${\rtheta} _-$ if and
only if $\epsilon (\theta \chi ^{-1},\psi_{ 0})=\epsilon
(\overline{\theta} \chi ^{-1},\psi_{ 0})=-1$, $\chi$ appears in
${\rtheta} ' _+$ if and only if\\ $\epsilon (\theta \chi ^{-1},\psi_{
0})=1$ and $\epsilon (\overline{\theta} \chi ^{-1},\psi_{ 0})=-1$,
and $\chi$ appears in ${\rtheta} ' _-$ if and only if $\epsilon
(\theta \chi ^{-1},\psi_{ 0})=-1$ and $\epsilon (\overline{\theta}
\chi ^{-1},\psi_{ 0})=1$.
\end{theorem}
Here $D^*_F$ is the unique quarternion division algebra over $F$.
This result was proved only in the odd residue characteristic case
in \cite{prasad}. Proof in the even residue characteristic case
appeared independently in \cite{prasad2} and \cite{vishnu}.

 We have, by definition of central
character, $\theta|_{F^*}=\omega_{{\rtheta}}\omega$ where
$\omega_{{\rtheta}}$ is the central character of ${\rtheta}$ and
$\omega=\omega_{K/F}$. Since $\theta|_{F^*}\neq \omega_{{\rtheta}}$
the character $\theta$ cannot occur in ${\rtheta}|_{K^*}$. A
necessary condition for a character $\lambda$ of $K^*$ to occur in
${\rtheta}|_{K^*}$ is that its restriction to $F^*$ should be equal
to the central character $\omega_{{\rtheta}}$. The question we would
like to ask at this point is whether $\theta$ twisted by some
character $\lambda$ of $K^*$ can occur in ${\rtheta}$ with
$\lambda|_{F^*}=\omega$. Note that such a twist satisfies the said necessary
condition.  Making it more precise, it means, whether there exist some $\lambda$ such that
$\lambda\theta$ occurs in ${\rtheta}|_{K^*}$. We prove some results
which give an affirmative answer to this question. In fact, we try
to count at each conductor level precisely how many characters occur
in ${\rtheta}_+$ and  ${\rtheta}_-$. It is not really surprising to
see that the necessary condition is not sufficient to guarantee the
occurrence of a character. Our computations on the local $\epsilon$-
factors are sometimes long, but by no means they are complicated. We feel that we have performed all kinds of computations possible using the $\epsilon$- factors of characters. Lending the words of Tunnel(\cite{tunnell}) the results here in this exposition are presented as an "entertainment". The main results in the exposition are coming in the last two sections.

\section{Notations}
Our notations are consistent with those used in \cite{vishnu} more
or less because we depend heavily on not only the results in
\cite{vishnu}, but also the computations performed there.

Throughout this paper $F$ will be a nonarchimedian local field of
characteristic $\not =2$ and $K$ a quadratic extension of $F$. The
image of $x\in K$ under the nontrivial element of the Galois group
of $K$ over $F$ is denoted by $\overline x$. For a local field $F$,
$O_F$ will be the ring of integers in $F$, $P_F=\pi_F O_F$ the
unique prime ideal in $ O_F$ and $\pi_F$ a uniformizer, i.e., an
element in $ P_F$ whose valuation is one, i.e., $v_F(\pi_F)=1$. The
cardinality of the residue field of $F$ is denoted by $q$ and $U_F=
O_F- P_F$ is the group of units in $ O_F$. Let $ P_F^i=\{x\in F
:v_F(x)\geq i\}$ and for $i\geq 0$ define $U_F^i= 1+ P_F^i$ (with
the proviso that $ U_F^0= U_F$).

Conductor of an additive character $\psi$ of $F$ or $K$ is $n(\psi)$
if $\psi$ is trivial on $P^{-n(\psi)}$, but nontrivial on
$P^{-n(\psi)-1}$. Fix an additive character $\psi$ of $F$ of
conductor zero (with no loss of generality, as in \cite{vishnu}) and
let $\psi_K= \psi \circ tr_{K/F}$ where $tr_{K/F}$ or simply $tr$ is
the trace map from $K$ to $F$. By $N_{K/F}$ or simply $N$ we mean
the norm map from $K$ to $F$ and by $d_{K/F}$ or simply $d$ the
differential exponent of $K$ over $F$ which is such that $tr
P_K^{-d}\subseteq O_F$ but $trP_K^{-d-1}\not \subseteq O_F$. The
conductor of $\psi_K$ is $d$. For a character $\chi$ of $F^*$ or
$K^*$ by $a(\chi )$ we mean the conductor of $\chi$, i.e., $a(\chi
)$ is the smallest integer $n\geq 0$ such that $\chi $ is trivial on
$ U^n$. We say that $\chi $ is unramified if $a(\chi )$ is zero.
Also, if $\chi_1$ and $\chi_2$ are two characters of $F$ then
$a(\chi_1\chi_2)\leq max(a(\chi_1),a(\chi_2))$. Equality holds if
$a(\chi_1)\neq a(\chi_2)$. Furthermore, $a(\chi)=a(\chi^{-1})$.

A character $\theta$ of $K^*$ is regular if it does not factor
through the norm map from $K$ to $F$. This guarantees that
$\theta\not=\overline{\theta}$.  The $F$-valuation of 2, $v_F(2)$,
will always be denoted by $t$. Therefore, $2=\pi_F^tu$, $u\in U_F$.
By $x_0$ we will always denote a nonzero element of $K$ with trace
$0$. Define $\psi_0$ by $\psi _0(x)=\psi(tr[-xx_0/2])$ for $x\in K$.
Then $\psi_0$ is an additive character of $K$ trivial on $F$.

If $G$ is a locally compact abelian group by $\widehat{G}$ we mean
the group of characters of $G$. Denote by $\omega_{K/F}$, or simply
$\omega$  the character of $F^*$ associated to $K$ by class field
theory, i.e., it is the unique nontrivial character of $F^*/N(K^*)$.

If $X$ is a finite set, by $|X|$ we will mean the number of elements
in $X$.

\section{Some useful results}
Deligne\cite{deligne} described how the epsilon factor changes under
twisting by a character of small conductor in the theorem:
\begin{theorem}\label{thdeligne}
Let $\alpha ,\beta $ be two characters of a local field $F$ such
that $a(\alpha )\geq 2a(\beta )$. Let $y_{\alpha}$ be an element of
$F^*$ such that $\alpha (1+x)=\psi (y_{\alpha}x)$ for $v_F(x)\geq
\frac{a(\alpha )}{2}$ (if $a(\alpha )=0$ let
$y_{\alpha}=\pi_F^{-n(\psi )}$). Then $\epsilon (\alpha \beta ,\psi
)=\beta^{-1}(y_{\alpha})\epsilon (\alpha ,\psi )$.
\end{theorem}
Note that $v_F(y_{\alpha})=-a(\alpha)-n(\psi)$.\\
From \cite{vishnu} we have
\begin{lemma}
If $a(\chi)\geq 2a(\ot)$, $\chi |_{F^*}=\omega$ then
\begin{equation}\label{epsdeligne}
\epsilon(\chi,\psi_0)=  \tilde{\omega}(-x_0/2)
\tilde{\omega}^{-1}(y_{\chi})
\end{equation}
where $\chi.\tilde{\omega}^{-1}(1+x)=\psi_K(y_{\chi{\oti}}x)$.

\end{lemma}
Here $y_{\chi}$ is as in theorem (\ref{thdeligne}).

The \emph{main theorem} in \cite{vishnu} states the following:

\begin{theorem}\label{mainresmain}
Let $K$ be a separable quadratic extension of a local field $F$ of
characteristic not two. Let $\psi$ be a nontrivial additive
character of $F$, and $x_0\in K^*$ such that $tr(x_0)=0$. Define an
additive
character $\psi_0$ of $K$ by $\psi_0(x)=\psi(tr[-xx_0/2])$. Then\\
\begin{equation}
\epsilon(\omega,\psi)\frac{\omega(\frac{x-\overline x}{x_0 -
\overline {x}_0})}{{{\left |{\frac{(x-\overline x)^2}{x\overline
x}}\right |}_{F^*}^{\frac{1}{2}}}} = \sum_{\chi \in S}\chi (x)
\end{equation}
$x\in K^*-F^*$ where as is usual, the summation on the right is by
partial sums over all characters of $K^*$ of conductor $\leq n$.
\end{theorem}

  We have the following result obtained by combining  corollary (7.2)  and the calculations given at the end of section 7 in \cite{vishnu}.

\begin{theorem}\label{sumclass}
Let $x=1+\pi_F^{r-1}\pi_Kx'$ where $x'\in U_F$, then
\[ \sum_{\chi\in
S(2r+2m)}\chi(x)= \left\{
\begin{array}{l l}
 -q^{r-1} & \quad \mbox{if $m=0$}\\
  0 & \quad \mbox{if $m=1,2,\ldots$}\\
   & \quad \mbox{and $m\neq d-1$}\\
\omega(-1)\epsilon(\omega,\psi)\frac{\omega(\frac{x-\overline x}{x_0
- \overline {x}_0})}{{{\left |{\frac{(x-\overline x)^2}{x\overline
x}}\right |}_{F^*}^{\frac{1}{2}}}} & \quad \mbox{if $m=d-1$}\\
\end{array} \right. \]

and
\[ \sum_{\chi '\in
S'(2r+2m)}\chi'(x)= \left\{
\begin{array}{l l}
 -q^{r-1} & \quad \mbox{if $m=0$}\\
  0 & \quad \mbox{if $m=1,2,\ldots$}\\
   & \quad \mbox{and $m\neq d-1$}\\
-\omega(-1)\epsilon(\omega,\psi)\frac{\omega(\frac{x-\overline
x}{x_0 - \overline {x}_0})}{{{\left |{\frac{(x-\overline
x)^2}{x\overline
x}}\right |}_{F^*}^{\frac{1}{2}}}} & \quad \mbox{if $m=d-1$}\\
\end{array} \right. \]

\end{theorem}

Note that \cite{vishnu} defined  $S$ to be the set $\{\chi \in K^* :
\chi \mid _{F^*} =\omega, \epsilon(\chi,\psi_0)=1\}$ and
$S(l)=\{\chi \in S:a(\chi)=l\}$. Analogously  defined $S'$ and
$S'(l)$ with the property that $\epsilon(\chi ,\psi_0)=-1$. For
computational convenience, we slightly changed our definition of $S$
to denote the set $=\{\chi \in \widehat{K^*} : \chi \mid _{F^*}
=\omega, \epsilon(\chi ^{-1},\psi_0)=1\}$ and $S(l)=\{\chi \in
S:a(\chi)=l\}$. Analogously, we redefined $S'$, and $S'(l)$. Because
of this change in notations, we have an extra term $\omega(-1)$ in
the above version compared to the one appeared in \cite{vishnu}.
This is due to the fact that  $\chi\overline{\chi}=1$ since their
restriction to $F^*$ is $\omega$ and so $\epsilon(\chi
^{-1},\psi_0)=\omega(-1)\epsilon(\chi ,\psi_0)$.
We define $S_l=S(l)\cup S'(l)$.\\
When $K/F$ is ramified, the following result can be verified
trivially by applying  lemma (5.1) in \cite{vishnu}.
\begin{lemma}\label{green}
Let $\chi \in \widehat{F^*}$ and $\psi$  a nontrivial character of
$(F,+)$.
\begin{enumerate}\label{guava}

\item{ If $n<a(\chi)+n(\psi)$ then
\begin{equation*}
\sum_{u\in\frac{U_F}{U_F^{a(\chi)}}}\chi^{-1}(u)\psi(\pi_F^{-n}u)=0
\end{equation*}}
\item{ If $n>a(\chi)+n(\psi)$ then
\begin{equation*}
\sum_{u\in\frac{U_F}{U_F^n}}\chi^{-1}(u)\psi(\pi_F^{-n}u)=0
\end{equation*}
}
\end{enumerate}
\end{lemma}

We also have the following theorem from \cite{vishnu}.
\begin{theorem}\label{indigo}
$|S(l)|=|S'(l)|$ for each feasible $l$, that is when $l=2d-1$ or
$l=2f$ with $f\geq d$.
\end{theorem}

We use theorem (\ref{prasad}) to determine whether a $\chi\in S$ is
such that $\chi\theta$ occurs in ${\rtheta}_+$ or ${\rtheta}_-$. By
this theorem, $\chi\theta$ occurs in ${\rtheta}_+$ if and only if
$\epsilon(\theta(\chi\theta)^{-1},\psi_0)=\epsilon(\chi^{-1},\psi_0)=1=\epsilon(\overline{\theta}(\chi\theta)^{-1},\psi_0)=
\epsilon(\chi^{-1}\tsbt,\psi_0)$ and $\chi\theta$ occurs in
${\rtheta}_-$ if and only if
$\epsilon(\chi^{-1},\psi_0)=-1=\epsilon(\chi^{-1}\tsbt,\psi_0)$.
Note that a character $\chi\theta$ can occur in $\rtheta$ if and
only if it occurs in either ${\rtheta}_+$ or in ${\rtheta}_-$. Also,
if $\chi \in S(l)$ for some $l$ then $\chi\theta$ can occur in
${\rtheta}$ if and only if it occurs in ${\rtheta}_+$. Furthermore
if $\chi\in S(l)$, then $\chi\theta$ cannot occur in ${\rtheta}_-$
since for that $\chi$, $\epsilon(\chi ^{-1},\psi_0)=+1$. Since its
multiplicity cannot exceed 1 in $\rtheta$, it is so in ${\rtheta}_+$
and ${\rtheta}_-$.

Now we are ready to start our counting. We divide the proof into
mainly two cases: $K/F$ ramified and $K/F$ unramified.
\section{Counting the twists when $K/F$ is ramified}

It is known ( see, for instance, \cite{tandon}, section 3) that if
$d$ is odd then $d=2t+1$ and there exists a uniformizer, denoted by
$\pi_K$ such that $tr\, \pi_K=0$. Let $x_0=\pi_K$. In this case
$\pi_K^2$ is a uniformizer of $F$ which we denote by $\pi_F$ and
$N\pi_K=-\pi_F$. If $d$ is even (which can only happen if the
residue characteristic is 2) then $O_K=O_F[\pi_K]$ where $\pi_K$ is
a uniformizer of $K$ which satisfies the Eisenstein polynomial
$X^2-u'\pi_F^sX-\pi_F$ with $s\leq t$. Again $N\pi_K=-\pi_F$. In
this case $d=2s$ and $\pi _K=\frac{\pi _F^su'}{2}(1+x_0)$ where
$x_0$ is a unit of trace $0$. We note that $n(\psi_0)$ is equal to 2
if $d$ is odd and $2(s-t)$ if $d$ is even. So $n(\psi_0)$ is always
even. Note also that if $\chi|_{F^*}=\omega$, then $a(\chi)$ is
either $2d-1$ or it is even, say $2f$, with $f\geq d$.

We know that (see \cite{vishnu}) if $\chi \in \widehat{F^*}$ and
$\psi$ is a nontrivial additive character of $F$, then

\begin{eqnarray}
\epsilon(\chi ,\psi) = \chi(c)q^{-a(\chi )/2}\sum _{y\in
\frac{U_F}{U_F^{a(\chi )}}}\chi ^{-1}(y)\psi(y/c)\label{mango}
\end{eqnarray}
where $v_F(c)=a(\chi)+n(\psi)$. In particular since $a(\omega)=d$
and we have chosen $\psi$ such that $n(\psi)=0$ we have
\begin{equation}
\label{banana} \epsilon(\omega , \psi)=\omega (\pi _F
^d)q^{-d/2}\sum_{y \in \frac{U_F}{U_F^d}}\omega (y)\psi(\pi _F
^{-d}y)
\end{equation}
This expression is obtained by normalizing the Haar measure given in
the expression for $\epsilon$-factor in \cite{tate} such that volume
of $O_F$ is 1.

To start off, we have the following simple lemma.
\begin{lemma}
For a regular character $\theta$ of $K^*$, $a(\frac{\theta}
{\overline{\theta}})$ is always even.
\end{lemma}
\begin{proof}
If not, suppose $a(\frac{\theta} {\overline{\theta}})=2r+1,\,r\geq
0$. Then it has to be nontrivial on $\frac{U_K^{2r}}{U_K^{2r+1}}$.
But $\frac{\theta}
{\overline{\theta}}(1+\pi_F^{r}a)=\frac{\theta(1+\pi_F^{r}a)}{\overline{\theta}(1+\pi_F^{r}a)}=1$,
where $a\in \mathbb{F}_q$ which is a contradiction.
\end{proof}

\subsection{Twist by characters of odd conductor}
We reserve the symbols $\ot$ and ${\ot}_{K/F}$ to denote elements of
$S_{2d-1}$.
\begin{lemma}\label{notwist}
If  $\frac{\theta} {\overline{\theta}}=(-1)^{v_K}$ then no
$\tilde{\omega}\theta$ can occur in ${\rtheta}$.
\end{lemma}
\begin{proof}
By theorem (\ref{prasad}), $\tilde{\omega}\theta$ can occur in
${r_{\theta}}_+$ if and only if $\epsilon (\tsbt\oti ,\psi_{\tiny
0})=\epsilon (\oti,\psi_{\tiny 0})=1$ and in ${r_{\theta}}_-$ if and
only if  $\epsilon (\tsbt\oti ,\psi_{\tiny 0})=\epsilon
(\oti,\psi_{\tiny 0})=-1$. Since $\tsbt$ unramified, we have
\begin{eqnarray*}
\epsilon (\tsbt\oti
,\psi_0)&=&\tsbt(\pi_K)^{a(\oti)+n(\psi_0)}\epsilon
(\oti,\psi_{\tiny
0})\\
&=&\tsbt(\pi_K)^{2d-1}\epsilon (\oti,\psi_{\tiny 0})\mbox{ (since $n(\psi_0)$ even)}\\
&=&-\epsilon (\oti,\psi_{\tiny 0})
\end{eqnarray*}
 which shows that
$\epsilon(\tilde{\omega}^{-1},\psi_0)=-\epsilon(\tilde{\omega}^{-1}\frac
{\overline{\theta}}{\theta},\psi_0)\,\forall\, \tilde{\omega}\in
S(2d-1)$. Similarly
$\epsilon(\tilde{\omega}^{-1},\psi_0)=-\epsilon(\tilde{\omega}^{-1}\frac
{\overline{\theta}}{\theta},\psi_0)$ for all $ \ot\in S'(2d-1)$. So
$\tilde{\omega}\theta$ can occur neither in ${r_\theta}_+$ nor in
${r_\theta}_-$ for any $\ot\in S_{2d-1}$. Therefore it cannot occur
in $\rtheta.$
\end{proof}

\begin{theorem}\label{mainoddconductor}
Let $0\neq a(\frac{\theta}{\overline{\theta}})<a(\tilde{\omega})$.
Then among  all $\tilde{\omega}\in S(2d-1)$ half and only half will
be such that $\tilde{\omega}\theta$ occur in ${r_{\theta}}_+$ and
among  all $\tilde{\omega}\in S'(2d-1)$ half and only half will be
such that $\tilde{\omega}\theta$ occur in ${r_{\theta}}_-$.
\end{theorem}
\rem When $d=1$, $a(\tilde{\omega})=1.$ Therefore, since
$a(\tbts)\neq 0$, this theorem is not applicable in $d=1$ case.
\begin{proof} We show that $\displaystyle{
\sum_{\tilde{\omega}\in S(2d-1)}\epsilon(\tilde{\omega}^{-1}\frac
{\overline{\theta}}{\theta},\psi_0)=0}$ so that half of $\ot\in
S(2d-1)$ will be such that $\etheta=+1$ and the other half $-1$. The
first half will occur in ${\rtheta}_+$. The remaining half will not
occur either in ${\rtheta}_+$ or in  ${\rtheta}_-$. The other part
of the proof is similar.\\

Note that $\nsi0$ is always even irrespective of $d$. Also,
$a(\oti)=a(\otitheta)$. Taking
$c=\pi_F^{d+\frac{\nsi0}{2}}\pi_K^{-1}$, in equation (\ref{mango})
we have if $\ot\in S(2d-1)$, then
\begin{eqnarray*}
\etheta&=&q^{-\frac{2d-1}{2}}\oti\tsbt(\pi_F^{d+\frac{\nsi0}{2}}\pi_K^{-1})\\
& &\times \sum_{y\in\frac{U_K}{U_K^{2d-1}}}\tilde{\omega}\tbts(y)\psi_0(\pi_F^{-(d+\frac{\nsi0}{2})}\pi_Ky)
\end{eqnarray*}
Write  $y\in\frac{U_K}{U_K^{2d-1}}$ as
$y=y_1(1+\pi_F^{r-1}\pi_Ky_2)$, $r\geq1,y_1\in \frac{U_F}{U_F^d}$,
$y_2=0$ or $y_2\in \frac{U_F}{U_F^d}$. Also note that
$\frac{\theta}{\overline {\theta}}$ is trivial on $F^*$. Summing
over $S(2d-1)$, we get
\begin{eqnarray*}
\lefteqn{\sum_{\ot\in
S(2d-1)}\etheta=}\\
& &q^{(-\frac{2d-1}{2})}\omega(\pi_F^{d+\frac{\nsi0}{2}})\tsbt(\pi_K)\sum_{y_1,y_2,r,\ot}[\ot(\pi_Ky_1(1+\pi_F^{r-1}\pi_Ky_2))\\
& &\tbts(1+\pi_F^{r-1}\pi_Ky_2)\psi_0(\pi_F^{-(d+\frac{\nsi0}{2})}\pi_Ky_1(1+\pi_F^{r-1}\pi_Ky_2))]
\end{eqnarray*}
But from the identity in theorem (\ref{mainresmain}) and and the
fact that we have to only consider characters in $S$ with odd
conductor (which is equal to $2d-1$) when $v_K(x)=1$, it follows
that


\begin{align*}
\lefteqn{\sum_{\ot}\ot(\pi_Ky_1(1+\pi_F^{r-1}\pi_Ky_2))= } & &\\
& &\left\{
\begin{array}{l l}\omega(-1)\epsilon(\omega,\psi)q^t \omega(y_1) \mbox{ \quad if $d=2t+1$}\\
  \omega(-1) \epsilon(\omega,\psi)q^{s-\frac{1}{2}}\times \omega(\pi_F^{s-t}uu'y_1(1+\pi_F^{s+r-1}u'y_2))  \mbox{ if $d=2s$}
  \end{array}
 \right .
 \end{align*}

 Also,
\begin{align*}
  \lefteqn{\psi_0(\pi_F^{-(d+\frac{\nsi0}{2})}\pi_Ky_1(1+\pi_F^{r-1}\pi_Ky_2)) = } & &\\
&  &\left\{
\begin{array}{l l}
\psi(-\pi_F^{-d}y_1) \mbox{ if $d=2t+1$}\\
\psi(-\pi_F^{-d}u^{-1}u'x_0^2y_1(1+\pi_F^{r+s-1}u'y_2))
\mbox{ if $d=2s$}\\
\end{array}
 \right.
 \end{align*}
Let $d=2t+1$.  If we keep $y_2$ fixed,
\begin{eqnarray*}
& &\sum_{y_1,\,\ot}\ot(\pi_Ky_1(1+\pi_F^{r-1}\pi_Ky_2))
\psi_0(\pi_F^{-(d+\frac{\nsi0}{2})}\pi_K(1+\pi_F^{r-1}\pi_Ky_2))\\
&=&\epsilon(\omega,\psi)q^t\sum
_{y_1}\omega(-y_1)\psi(-\pi_F^{-d}y_1)\\
&=&\epsilon(\omega,\psi)q^t \epsilon(\omega,\psi)\omega(\pi_F^d)q^d
\end{eqnarray*}
which is a multiple of $\epsilon(\omega,\psi)$ independent of $y_1$
and $y_2$. So
\begin{eqnarray*}
\sum_{\ot\in
S(2d-1)}\etheta=C\sum_{r,y_2}\tbts((1+\pi_F^{r-1}\pi_Ky_2) =C\times
0=0
\end{eqnarray*}
since $\tbts$ is a nontrivial character of
$\frac{U_K}{U_FU_K^{a(\tbts)}}$ and $a(\tbts)\leq 2d-2$. Here $C$ is
a constant independent of $y_1$ and $y_2$. Similarly if $d=2s$, then
if we again keep $y_2$ fixed

 \begin{align*}
\lefteqn{\sum_{y_1,\,\ot}\ot(\pi_Ky_1(1+\pi_F^{r-1}\pi_Ky_2))\psi_0(\pi_F^{-(d+\frac{\nsi0}{2})}\pi_K(1+\pi_F^{r-1}\pi_Ky_2)y_1)}\\
\lefteqn{=\epsilon(\omega,\psi)q^{s-\frac{1}{2}}\omega(-\pi_F^{s-t}uu')\times}& &\\
& &\sum_{y_1}\omega((1+\pi_F^{s+r-1}u'y_2)y_1)\psi(-\pi_F^{-d}u^{-1}u'x_0^2y_1(1+\pi_F^{r+s-1}u'y_2))
 \end{align*}
which is again a constant multiple of $\epsilon(\omega,\psi)$
independent of $y_1$ and $y_2$. So
\begin{eqnarray*}
\sum_{\ot\in
S(2d-1)}\etheta=C'\sum_{r,y_2}\tbts((1+\pi_F^{r-1}\pi_Ky_2)
=C'\times 0=0
\end{eqnarray*}
where $C'$ is a constant multiple of $\epsilon(\omega,\psi)$. This
completes the proof of the theorem.
\end{proof}

\begin{coro}\label{odd1}
The number of $\ot\in S_{2d-1}$ such that $\tilde{\omega}\theta$
occurs in $r_{\theta}$ is $|S_{2d-1}|/2=|S(2d-1)|=|S'(2d-1)|$.
\end{coro}

\begin{proof}
This is clear since occurring in $\rtheta$ means occurring in either
${\rtheta}_+$ or in ${\rtheta}_-$. Equality follows from theorem
(\ref{indigo}).
\end{proof}
\begin{lemma}\label{odd2}
If $a(\tbts)>a(\ot)$ then the number of $\ot\in S_{2d-1}$ such that
$\tilde{\omega}\theta$ occurs in $r_{\theta}$ is $|S_{2d-1}|/2$.
\end{lemma}
\begin{proof}
This is quite easy to verify.  In this case,
$a(\tbts)=a(\oti\tbts)$. Note that $a(\tbts)$ is even. So if
\begin{equation}\label{notequal}\eoti\neq\etheta,
\end{equation}
 consider the character $\mu=(-1)^{v_K}$ of $K^*$ and take
$\oti_2=\oti\mu$. If we consider the expression for epsilon factors
on both the sides of (\ref{notequal}), since $a(\oti\tsbt)$ is even,
no $\pi_K$ is present but only $\pi_F$ on the RHS of this equation.
Therefore the twist by $\mu$ will not make any difference on the
RHS. But on the LHS, an extra $\mu(\pi_K)=-1$ will appear changing
the sign of LHS. Similarly if $\eoti=\etheta$, we can make them
unequal by the same sort of twisting. So for half of $\ot\in
S_{2d-1}$, the corresponding epsilon factors are equal and for the
other half they are unequal.
\end{proof}

\subsection{Twist by characters of even conductor}
Note that if $a(\lambda)=2f\geq 2d$, then in the expression for
$\epsilon(\lambda ^{-1},\psi_0)$ there is no $\pi_K$, but only
$\pi_F$.
\begin{theorem}

Let $\lambda\in S(2f+2d),\, f\geq 0\quad a(\tbts)\leq
a(\lambda)-2d=2f$. Then all the elements in $\{\lambda\theta:
\lambda \in S(2f+2d)\}$ will occur in ${\rtheta}_+$. Similarly if
$\lambda'\in S'(2f+2d)$, then all the elements in $\{\lambda'\theta:
\lambda' \in S'(2f+2d)\}$ will occur in ${\rtheta}_-$. Therefore the
number of $\lambda\theta$ where $\lambda \in S_{2f+2d}$ occurring in
$\rtheta$ is $|S_{2f+2d}|$.
\end{theorem}
\begin{proof}
Consider the two sums \\$\displaystyle {\sum_{\lambda\in
S(2f+2d)}\eli}$ and $\displaystyle {\sum_{\lambda\in
S(2f+2d)}\eltheta},$ \\

we have
\begin{align*}
\lefteqn{\sum_{\lambda\in S(2f+2d)}\eltheta=}& &\\
& &q^{-f-d}\omega(\pi_F^{f+d+\frac{\nsi0}{2}})\sum_{\lambda\in
S(2f+2d)}\sum_{y\in\frac{U_K}{U_K^{2f+2d}}}\lambda(y)\tbts(y)\psi_0(\pi_F^{-(f+d+\frac{\nsi0}{2})}y)
 \end{align*}
In this summation, by theorem (\ref{sumclass}),
$\displaystyle{\sum_{\lambda\in S(2f+2d)}\lambda(y)\neq0}$ only for
two types of $y$'s: \\
1). when $y=y_1(1+\pi_F^f\pi_Ky_2)$ where
$y_1,y_2\in U_F$, and \\
2). when $y=y_1(1+\pi_F^{f+d-1}\pi_Ky_2)$
where $y_1,y_2\in U_F$ or $y_2=0$. \\
But since $a(\tbts)\leq 2f$ and
$\tbts=1$ on $F^*$ we have $\tbts$ trivial on these $y$'s. So both
the sums are independent of $\tbts$ and so they are the same. That
is,
$$\sum_{\lambda\in S(2f+2d)}\eli =\sum_{\lambda\in
S(2f+2d)}\eltheta$$ which means $$\eli =\eltheta\quad \forall
\lambda \in S(2f+2d)$$ since $\eli=1$ for each $\lambda\in
S(2f+2d)$. The remaining part follows similarly.
\end{proof}
\begin{coro}
If $\tbts=(-1)^{v_ K}$, then all $\lambda\in S(2f+2d)$  are such
that $\lambda\theta$  occur in ${\rtheta}_+$. Similarly all
$\lambda'\in S'(2f+2d)$  are such that $\lambda'\theta$  occur in
${\rtheta}_-$.
\end{coro}
\begin{proof}
It follows by taking $a(\tbts)=0$ in the above theorem.
\end{proof}
Note: The above corollary shows the differrence between characters
of even conductor and characters of odd conductor. This corollary is
extremely opposite to  lemma (\ref{notwist}).\\
Let $\lambda\in S(2f+2d),\,2f<a(\tbts)<a(\lambda)$. Note that if
$d=1$, then no such $\theta$ exists. So we have $d\geq 2$ and so $q$
even. By definition, we have
\begin{align*}
\lefteqn{\eltheta=}& &\\
& & q^{-f-d}\omega(\pi_F^{f+d+\frac{\nsi0}{2}})\sum_{\lambda\in
S(2f+2d)}\sum_{y\in\frac{U_K}{U_K^{2f+2d}}}\lambda(y)\tbts(y)\psi_0(\pi_F^{-(f+d+\frac{\nsi0}{2})}y)
\end{align*}
Again, by theorem (\ref{sumclass}), the sum $\sum_{\lambda\in
S(2f+2d)}\lambda(y)\neq 0$ only for two types of $y$'s:
\begin{enumerate}

\item{$y=y_1(1+\pi_F^{f+d-1}\pi_Ky_2), \, y_1\in \frac{U_F}{U_F^{f+d}},\, y_2 \in \mathbb{F}_q$}
\item{$y=y_1(1+\pi_F^{f}\pi_Ky_2), \, y_1\in \frac{U_F}{U_F^{f+d}},\, y_2 \in \frac{U_F}{U_F^d}$}
\end{enumerate}
Consider first type of $y$'s. $\tbts$ is trivial on the these $y$'s.
Now $$\sum_{\lambda\in
S(2f+2d)}\lambda(y)=\omega(y_1)\sum_{\lambda\in
S(2f+2d)}\lambda(1+\pi_F^{f+d-1}\pi_Ky_2)=-q^{f+d-1}\omega(y_1)$$

Also
 \begin{align*}
\psi_0(\pi_F^{-f-d-\frac{\nsi0}{2}}y)= \left\{
\begin{array}{l l}\psi(-\pi_F^{-1}y_1y_2) \mbox{ \quad if $d=2t+1$}\\
  \psi(-\pi_F^{-1}u^{-1}u'x_0^2y_1y_2) \mbox{ if $d=2s$}
  \end{array}
 \right .
 \end{align*}
Therefore
\begin{eqnarray*}
&&\sum_{\lambda\in
S(2f+2d)}\lambda(y)\psi_0(\pi_F^{-f-d-\frac{\nsi0}{2}}y)\\
&&=\begin{cases}-q^{f+d-1}\sum_{y_2\in\mathbb{F}_q}\sum_{y_1\in\frac{U_F}{U_F^{f+d}}}\omega(y_1
)\psi(-\pi_F^{-1}y_1y_2)\\
\qquad\mbox{ if $d$ odd} \cr
-q^{f+d-1}\sum_{y_2\in\mathbb{F}_q}\sum_{y_1\in\frac{U_F}{U_F^{f+d}}}\omega(y_1)\psi(-\pi_F^{-1}u^{-1}u'x_0^2y_1y_2)\\
\qquad\mbox{
if $d$ even}\end{cases}\\
&&=\begin{cases}-q^{f+d-1}\sum_{y_2\in\mathbb{F}_q}\omega(-y_2)\sum_{y_1\in\frac{U_F}{U_F^{f+d}}}\omega(y_1)\psi(\pi_F^{-1}y_1)\\
\qquad\mbox{
if $d$ odd} \cr
-q^{f+d-1}\sum_{y_2\in\mathbb{F}_q}\omega(-y_2x_0^2uu')\sum_{y_1\in\frac{U_F}{U_F^{f+d}}}\omega(y_1)\psi(\pi_F^{-1}y_1)\\
\qquad\mbox{
if $d$ even}\end{cases}
\end{eqnarray*}
Since $a(\omega)=d\neq 1$, by lemma ({\ref{guava}})
$\displaystyle{\sum_{y_1\in\frac{U_F}{U_F^{f+d}}}\omega(y_1)\psi(\pi_F^{-1}y_1)=0.}$
So $$\sum_{\lambda\in
S(2f+2d)}\lambda(y)\sum_{y_1\in\frac{U_F}{U_F^{f+d}}}\lambda(y)\psi_0(\pi_F^{-f-d-\frac{\nsi0}{2}}y)=0$$

Consider the second type of $y$'s. On these, we have
 \begin{align*}
\lefteqn{ \sum_{\lambda\in
S(2f+2d)}\lambda(y_1(1+\pi_F^{f}\pi_Ky_2))=}\\
&&\begin{cases}\omega(-1)\omega(y_1)\omega(\pi_F^fy_2)q^{f+t+\frac{1}{2}}\epsilon(\omega,\psi)\mbox{
if $d$ odd}\cr
\omega(-1)\omega(y_1y_2uu')\omega(\pi_F^{f+s-t}y_2)q^{f+s}\epsilon(\omega,\psi)\mbox{
if $d$ even}\end{cases}\\
&&\mbox{ by theorem (\ref{mainresmain}) and theorem
(\ref{sumclass})}
\end{align*}
and
\begin{eqnarray*}
\psi_0(\pi_F^{-f-d-\frac{n(\psi_0)}{2}}y_1(1+\pi_F^{f}\pi_Ky_2))=\begin{cases}\psi(-\pi_F^{-d}y_1
y_2)\mbox{ if $d$ odd}\cr \psi(-\pi_F^{-d}y_1
y_2u^{-1}u'x_0^2)\mbox{ if $d$ even}\end{cases}
\end{eqnarray*}
Let $d=2t+1$. Then
\begin{eqnarray*}
&&\sum_{\lambda\in S(2f+2d)}\sum_{y_1 \in
\frac{U_F}{U_F^{f+d}}}\sum_{y_2 \in
\frac{U_F}{U_F^d}}\lambda(y)\tbts(y)\psi_0(\pi_F^{-f-d-1}y)\\
&&=\omega(-1)\omega(\pi_F^f)q^{f+t+\frac{1}{2}}\epsilon(\omega,\psi)\\
&&\qquad\times\sum_{y_1
\in \frac{U_F}{U_F^{f+d}}}\sum_{y_2 \in
\frac{U_F}{U_F^d}}\omega(y_1y_2)\psi(-\pi_F^{-d}y_1
y_2)\tbts(1+\pi_F^{f}\pi_Ky_2)\\
&&=q^{f}q^{f+t+\frac{1}{2}}\omega(\pi_F^f)\epsilon(\omega,\psi)\\
&&\qquad\times\sum_{y_2
\in \frac{U_F}{U_F^d}}\tbts(1+\pi_F^{f}\pi_Ky_2)\sum_{y_1 \in
\frac{U_F}{U_F^d}}\omega(y_1y_2)\psi(\pi_F^{-d}y_1 y_2)\\
&&=q^{2f}q^{\frac{d}{2}}\omega(\pi_F^f)\epsilon(\omega,\psi)\epsilon(\omega,\psi)\omega(\pi_F^d)q^{\frac{d}{2}}\sum_{y_2
\in \frac{U_F}{U_F^d}}\tbts(1+\pi_F^{f}\pi_Ky_2)\\
&&=q^{2f+d}\omega(\pi_F^{f+d})\epsilon(\omega,\psi)^2\sum_{y_2 \in
\frac{U_F}{U_F^d}}\tbts(1+\pi_F^{f}\pi_Ky_2)
\end{eqnarray*}
If $d=2s$ we will get the same sum with an extra $\omega(-1)$
factor. Now if $a(\tbts)\geq 2f+4$ then $$\sum_{y_2 \in
\frac{U_F}{U_F^d}}\tbts(1+\pi_F^{f}\pi_Ky_2)=0 \mbox { and so }
\sum_{\lambda \in S(2f+2d)}\eltheta=0.$$ Therefore half of elements
in $\{\lambda\theta:\lambda \in S(2f+2d)\}$ will appear in
${\rtheta}_+$. Similarly, half of
elements in $\{\lambda'\theta:\lambda '\in S'(2f+2d)\}$ will appear in ${\rtheta}_-$. \\
Let $a(\tbts)=2f+2$. Then \\
$\displaystyle{\sum_{y_2 \in
\frac{U_F}{U_F^d}}\tbts(1+\pi_F^{f}\pi_Ky_2)=q^{d-1}\sum_{a\in\mathbb{F}_q}\tbts(1+\pi_F^{f}\pi_K
a)=-q^{d-1}}.$

Therefore  if $d=2t+1$, then
\begin{align*}
\sum_{\lambda\in S(2f+2d)}\sum_{y_1 \in \frac{U_F}{U_F^{f+d}}}\sum_{y_2 \in \frac{U_F}{U_F^d}}\lambda(y)\tbts(y)\psi_0(\pi_F^{-f-d-1}y)&\\
=-q^{2f+2d-1}\omega(\pi_F^{f+d})&\epsilon(\omega,\psi)^2
 \end{align*}

So \begin{eqnarray*}
\sum_{\lambda\in
S(2f+2d)}\eltheta&=&q^{-f-d}\omega(\pi_F^{f+d})\times
-q^{2f+2d-1}\omega(\pi_F^{f+d})\epsilon(\omega,\psi)^2\\&=&-q^{f+d-1}\epsilon(\omega,\psi)^2
 \end{eqnarray*}

Similarly we will get $\displaystyle{\sum_{\lambda\in
S(2f+2d)}\eli=(q-1)q^{f+d-1}\epsilon(\omega,\psi)^2}$

(This is
because in place of $\displaystyle{\sum_{y_2 \in
\frac{U_F}{U_F^d}}\tbts(1+\pi_F^{f}\pi_Ky_2)}$ we have
$|\frac{U_F}{U_F^d}|$).

But $$\sum_{\lambda\in
S(2f+2d)}\eli=|S(2f+2d)|=(q-1)q^{f+d-1}.$$ So
$\epsilon(\omega,\psi)^2=1$. (If $d=2s$, instead of this, we have
$\omega(-1)\epsilon(\omega,\psi)^2=1$). Therefore number of
$\lambda$ such that $\lambda\theta$ appear in ${\rtheta}_+$ is

$\sum_{\lambda \in S(2f+2d)}(\eli +\eltheta)=
\frac{(q-1)-1}{2}q^{f+d-1}=\frac{q-2}{2}q^{f+d-1}$. If $d=2s$ then
also  we can show that the sum is $\frac{q-2}{2}q^{f+d-1}$.

So in ${\rtheta}_+$, the number of $\lambda\theta$ occurring  where
$\lambda\in S(2f+2d)$ is $\frac{q-2}{2}q^{f+d-1}$. Similarly in
${\rtheta}_-$, the number of $\lambda'\theta$ occurring  where
$\lambda'\in S'(2f+2d)$ is $\frac{q-2}{2}q^{f+d-1}$. We summarize
the above computations in the following two theorems.
\begin{theorem}
Let $\lambda\in S(2f+2d),\,2f+2<a(\tbts)<a(\lambda)$. Then among all
$\lambda\theta$ where $\lambda \in S(2f+2d)$ exactly half will occur
in ${\rtheta}_+$. Similarly, let $\lambda'\in
S'(2f+2d),\,2f+2<a(\tbts)<a(\lambda')$. Then among all
$\lambda'\theta$ where $\lambda' \in S'(2f+2d)$ exactly half will
occur in ${\rtheta}_-$. Therefore the number of $\lambda\theta $
where $\lambda \in S_{2f+2d}$ occurring in $\rtheta$ is
$|S_{2f+2d}|/2$.
\end{theorem}
\begin{theorem}
Let $\lambda\in S(2f+2d),\,a(\tbts)=2f+2<a(\lambda)$. Then number of
$\lambda\theta$ appearing in ${\rtheta}_+$ where $\lambda \in
S(2f+2d)$ is $\frac{q-2}{2}q^{f+d-1}$. Similarly, let $\lambda'\in
S'(2f+2d),\,a(\tbts)=2f+2<a(\lambda')$. Then number of
$\lambda'\theta$ appearing in ${\rtheta}_-$ where $\lambda' \in
S'(2f+2d)$ is $\frac{q-2}{2}q^{f+d-1}$. The number of $\lambda\theta
$ where $\lambda \in S_{2f+2d}$ occurring in $\rtheta$ is therefore
$(q-2)q^{f+d-1}$.
\end{theorem}
Note: These two theorems are not valid for $d=1$ since no $\tbts$
satisfies the condition in the theorem.
\begin{theorem}\label{even2}
Let $a(\lambda)=2f+2d<a(\tbts)=2m<a(\lambda)+2d$. Then the number of
$\lambda\theta$ with $\lambda\in S_{2f+2d}$ appearing in $\rtheta$
is $|S(2f+2d)|=|S_{2f+2d}|/2$.
\end{theorem}

\begin{proof}
Here $a(\lambda ^{-1}\tsbt)=a(\tbts)$. Using the definition of
$\epsilon$-factors, we have
\begin{eqnarray*}
\lefteqn{\sum_{\lambda\in S(2f+2d)}\eltheta=}& &\\
& &q^{-m}\omega(\pi_F^{m+\frac{\nsi0}{2}})\sum_{\lambda\in
S(2f+2d)}\sum_{y\in\frac{U_K}{U_K^{2m}}}\lambda(y)\tbts(y)\psi_0(\pi_F^{-(m+\frac{\nsi0}{2})}y)
\end{eqnarray*}
Recall that, from theorem (\ref{sumclass}) the sum $\sum_{\lambda\in
S(2f+2d)}\lambda(y)\neq 0$ only for three types of $y$'s:
\begin{enumerate}
\item{$y=y_1(1+\pi_F^{f}\pi_Ky_2), \, y_1\in \frac{U_F}{U_F^m},\, y_2 \in \frac{U_F}{U_F^{m-f}}$}
\item{$y=y_1(1+\pi_F^{f+d-1}\pi_Ky_2), \, y_1\in \frac{U_F}{U_F^m},\, y_2 \in \frac{U_F}{U_F^{m-f-d+1}}$}
\item{$y=y_1(1+\pi_F^{f+d}\pi_Ky_2), \, y_1\in \frac{U_F}{U_F^m},\, y_2 \in \frac{U_F}{U_F^{m-f-d}}$ or $y_2=0$}
\end{enumerate}
But on the third type of $y$'s, $\lambda$ is just $\omega$ since
$a(\lambda)=2f+2d$. On the second type of $y$'s,
\begin{eqnarray*}\sum_{\lambda\in S(2f+2d)}\lambda(y) &=& \omega(y_1)\sum_{\lambda\in
S(2f+2d)}\lambda(1+\pi_F^{f+d-1}\pi_Ky_2)\\
 &=&  \omega(y_1)(-q^{f+d-1}) \mbox{ (by theorem (\ref{sumclass})) }
\end{eqnarray*}
So $\displaystyle{\sum_{\lambda\in S(2f+2d)}\lambda(y)}$ is
independent of $\lambda$ on these $y$'s. Finally consider the first
type of $y$'s. Let $d=2t+1$.
\begin{eqnarray*}
\sum_{\lambda\in S(2f+2d)}\lambda(y_1(1+\pi_F^f \pi_Ky_2))&=&
\omega(y_1)\sum_{\lambda\in S(2f+2d)}\lambda(1+\pi_F^f \pi_Ky_2)\\
&=&\omega(-1)\omega(y_1) \epsilon(\omega,\psi)\omega(\pi _F^f
y_2)q^{f+t+\frac{1}{2}}\\
& & \mbox{by theorem (\ref{mainresmain}) and
theorem(\ref{sumclass}).}
\end{eqnarray*}
Therefore
\begin{align*}
\lefteqn{\sum_{\lambda\in S(2f+2d)}\sum_{y_1 \in \frac{U_F}{U_F^m}}\sum_{y_2 \in \frac{U_F}{U_F^{m-f}}}\lambda(y_1(1+\pi_F^f
\pi_Ky_2))\tbts(y_1(1+\pi_F^f \pi_Ky_2))}& &\\
& &\times \psi_0(\pi_F^{-m-1}y_1(1+\pi_F^f \pi_Ky_2)) \\
\lefteqn{=\epsilon(\omega,\psi)\omega(\pi _F^f
)q^{f+t+\frac{1}{2}}}& &\\
& &\times \sum_{y_1 \in \frac{U_F}{U_F^m}}\sum_{y_2 \in
\frac{U_F}{U_F^{m-f}}} \omega(-y_1 y_2)\psi(-\pi_F^{-m+f}y_1
y_2)\tbts (1+\pi_F^f \pi_Ky_2)
\end{align*}
In this sum, we have
\begin{align*}
\lefteqn{\sum_{y_1 \in \frac{U_F}{U_F^m}}\omega(y_1 y_2)\psi(-\pi_F^{-m+f}y_1
y_2)}& &\\
& &=\sum_{y_1 \in \frac{U_F}{U_F^{m-f}}}\omega(y_1
y_2)\psi(-\pi_F^{-m+f}y_1 y_2)|\frac{U_F^{m-f}}{U_F^m}|
\end{align*}

Since $a(\omega)=d$ and $m-f>d$, by lemma (\ref{guava}), the above
sum is zero.


Now if $d=2s$, we have $n(\psi_0)=2(s-t)$. Also, here the trace 0
element $x_0$ is a unit and  $\pi _K=\frac{\pi _F^su'}{2}(1+x_0)$.
Considering  $y$'s  first type, we have
\begin{eqnarray*}
\psi_0(\pi_F^{-(m+\frac{\nsi0}{2})}y)&=&
\psi_0(\pi_F^{-m-s+t}y_1(1+\pi_F^f \pi_Ky_2))\\
&=& \psi (-\frac{\pi_F^{-m-s+t}}{2} y_1 tr\,x_0(1+\pi_F^f
\pi_Ky_2))\\
&=& \psi (-\frac{\pi_F^{-m-s+t}}{2} y_1x_0^2.2\pi_F^f y_2\frac{\pi_F^s}{2}u')\\
&=&\psi(-\pi_F^{-m+f}u^{-1}u'y_1y_2x_0^2)
\end{eqnarray*}
and so
\begin{align*}
\lefteqn{\sum_{\lambda\in S(2f+2d)}\lambda(y_1(1+\pi_F^f \pi_Ky_2))=}& &\\
& &\omega(-1)\epsilon(\omega,\psi)\omega(\pi _F
^{f+s-t})q^{f+s}\omega(y_1 y_2uu')
\end{align*}
by theorem (\ref{mainresmain}) and
theorem(\ref{sumclass}).

Therefore
\begin{align*}
\lefteqn{\sum_{\lambda\in S(2f+2d)}\sum_{y_1 \in \frac{U_F}{U_F^m}}\sum_{y_2
\in \frac{U_F}{U_F^{m-f}}}[\lambda(y_1(1+\pi_F^f
\pi_Ky_2))\tbts(y_1(1+\pi_F^f \pi_Ky_2))\times}& &\\
& &\psi_0(\pi_F^{-m-s+t}y_1(1+\pi_F^f \pi_Ky_2))]\\
=\lefteqn{\epsilon(\omega,\psi)\omega(\pi _F ^{f+s-t})q^{f+s}\times}& & \\
& &\sum_{y_2 \in \frac{U_F}{U_F^{m-f}}} \tbts(1+\pi_F^f \pi_Ky_2)\sum_{y_1 \in \frac{U_F}{U_F^m}}\omega(y_1
y_2u^{-1}u')\psi(-\pi_F^{-m+f}u^{-1}u'y_1y_2x_0^2)
\end{align*}

The sum $\displaystyle{ \sum_{y_1 \in \frac{U_F}{U_F^m}}\omega(y_1
y_2u^{-1}u')\psi(\pi_F^{-m+f}u^{-1}u'y_1y_2)= 0}$ as in the $d=2t+1$
case since $m-f>d$. So the sum over the first type of $y$'s become
zero. So in both $d$ odd and even cases, the sum
$\displaystyle{\sum_{\lambda\in S(2f+2d)}\eltheta}$  depends only on
second and third type of $y$'s and is independent of $\lambda$.
Suppose this sum is $n$. Using similar arguments,  we have
$\displaystyle{\sum_{\lambda'\in S'(2f+2d)}\elitheta=n}.$ So the
number of $+1$'s in $\{\eltheta: \lambda\in
S(2f+2d)\}=\frac{|S(2f+2d)|+n}{2}$. Similarly, number of $-1$'s in
$\{\elitheta: \lambda'\in S'(2f+2d)\}=-\frac{-|S(2f+2d)|+n}{2}$.
Therefore, the number of $\lambda\theta$ appearing in ${\rtheta}_+$
is $\frac{|S(2f+2d)|+n}{2}$ , number of $\lambda'\theta$ appearing
in ${\rtheta}_-$ is $\frac{|S(2f+2d)|-n}{2}$. Total number of
$\lambda\theta$ appearing in $\rtheta$ is $|S(2f+2d)|$.
\end{proof}
When $a(\lambda)$ is too small compared to $a(\tbts)$  the
occurrence of $\lambda\theta$ in ${\rtheta} _+$  or ${\rtheta} _-$
depends only on $\theta$.
\begin{theorem}\label{even1}
Suppose $\lambda\in S(2m)$, $m\geq d$ and $a(\tbts)=2n \geq
a(\lambda)+2d$. Then either all the elements in $\{\lambda\theta :
\lambda \in S(2m)\}$ will occur in ${\rtheta} _+$ or all the
elements in $\{\lambda'\theta : \lambda' \in S'(2m)\}$ will occur in
${\rtheta} _-$ not both. Therefore the number of $\lambda\theta$
where $\lambda \in S_{2m}$ occurring in $\rtheta$ is $|S_{2m}|/2$.
\end{theorem}
\begin{proof}
We have if $\chi\in \widehat{K^*}$ with $\chi|_{K^*}=\omega$ and
$a(\chi)\geq 2a(\ot)$ then $ \epsilon(\chi,\psi_0) =
\tilde{\omega}(-x_0/2) \tilde{\omega}^{-1}(y_{\chi})$ where\\
$y_{\chi}=\begin{cases}\pi_F^{-f-\frac{d-1}{2}}\pi_Ka_0(\chi)(1+a_1(\chi)\pi_K)(1+a_2(\chi)\pi_F)...
& {\rm if}\  d\  {\rm is\  odd} \cr
\pi_F^{-f-\frac{d}{2}}x_0a_0(\chi)(1+a_1(\chi)\pi_K)(1+a_2(\chi)\pi_F)...
& {\rm if}\ d\  {\rm is\ even}\end{cases}$.

Here $a(\lambda ^{-1}\tsbt)=a(\tsbt)\geq 4d  >2a(\ot)$. Therefore

$\epsilon(\lambda^{-1}\tsbt,\psi_0)=\tilde{\omega}(-x_0/2)\times$\\
$\tilde{\omega}^{-1}(\pi_F^{-f-\frac{d-1}{2}}
\pi_Ka_0(\lambda^{-1}\tsbt)(1+a_1(\lambda^{-1}\tsbt)\pi_K)...(1+a_{2d-2}
(\lambda^{-1}\tsbt)\pi_F^{d-1}))$

if $d$ odd.

But note that
$(\lambda ^{-1}\tsbt)|_{U_K^{2n-2d+1}}$ determines $a_i(\lambda
^{-1}\tsbt)$ for $i=0,1,\ldots 2d-2$ and on $U_K^{2n-2d+1}$,
$\lambda ^{-1}\tsbt=\tsbt$. Therefore
$\epsilon(\lambda^{-1}\tsbt,\psi_0)$ is independent of $\lambda$ or
$\epsilon(\lambda^{-1}\tsbt,\psi_0)=\epsilon(\tsbt,\psi_0)$ Now
suppose that $\epsilon(\lambda^{-1},\psi_0)\neq
\epsilon(\tsbt,\psi_0)=-1$ for one $\lambda \in S(2m)$. Then for all
$\lambda' \in S'(2m),\,
\epsilon(\lambda'^{-1},\psi_0)=-1=\epsilon(\tsbt,\psi_0)=\epsilon(\lambda'^{-1}\tsbt,\psi_0)$.
Therefore $\{\lambda'\theta : \lambda' \in S'(2m)\}$ will occur in
${\rtheta} _-$. On the other hand if $\eli=\etsbt=1$ for one
$\lambda$  it is the same for all other $\lambda \in S(2m)$. This
proves the theorem.
\end{proof}

\begin{coro}
Suppose $a(\lambda)=2f+2d<a(\tbts)=2m$. If $n$= number of
$\lambda\theta$, $\lambda \in S(2f+2d)$, appearing in ${\rtheta}_+$
then number of $\lambda'\theta$, $\lambda' \in S'(2f+2d)$, appearing
in ${\rtheta}_-$ is $|S(2f+2d)|-n=|S'(2f+2d)|-n$. Also, if
$2m>a(\lambda)+2d$, then either $n=0$ or $n=|S(2f+2d)|$.
\end{coro}
\begin{proof} Follows easily from the above two theorems.
\end{proof}

 Only one case is left now for us to handle in this exposition viz.
$a(\tbts)=a(\lambda)$. In this case we are not giving an exact
count, but still we provide a lower bound in the next theorem. Note
that our calculations deal much with $a(\tbts\lambda)$ and it is
difficult to find when the two characters have equal conductor.
\begin{theorem}
If $a(\tbts)=a(\lambda)=2f+2d$, $\lambda|_{F^*}=\omega$ then the
number of $\lambda \theta$ appearing in $\rtheta$ is greater than or
equal to $q^{f+d-1}$.
\end{theorem}
\begin{proof}
Note that $S(2f+2d)\cup S'(2f+2d)=\{\tsbt\chi :
\chi|_{F^*}=\omega,a(\chi)=2d-1,2d,2d+2,\ldots, 2f+2d-2\}\bigcup
\{\tsbt\chi :
\chi|_{F^*}=\omega,a(\chi)=2f+2d,\,\chi|_{U_K^{2f+2d-1}}\neq
\tbts|_{U_K^{2f+2d-1}} \}$. Now a
$\tsbt\chi.\theta=\chi\overline{\theta}$ will appear in $\rtheta$ if
and only if
${\epsilon(\chi^{-1},\psi_0)}=\epsilon(\chi^{-1}\tbts,\psi_0)$. So
the number of $\chi\overline{\theta}$ appearing in $\rtheta$ where
$a(\chi)=2d-1,2d,\ldots, 2f+2d$ is greater than or equal to
$|S(2d-1)|+|S(2d)|+S(2d+2)|+\ldots+|S(2f+2d-2)|=q^{d-1}+(q-1)q^{d-1}+(q-1)q^d+\ldots+(q-1)q^{f+d-2}=q^{f+d-1}$
by corollary (\ref{odd1}), lemma (\ref{odd2}), and theorems
(\ref{even2}) and (\ref{even1}). Note that we are not considering
$\chi$'s with conductor equal to $2f+2d$ and that is why we are
unable to claim equality.
\end{proof}
\rem If $q=2$, there is no $\chi$ such that
$\chi|_{U_K^{2f+2d-1}}\neq \tsbt|_{U_K^{2f+2d-1}}$. So equality
holds in the theorem.

\section{The unramified case}
Suppose $K$ over $F$ is unramified and let $\chi\in \widehat{F^*}$
be such that  $\chi|_{K^*}=\omega$. Let $\ot$ be an extension of
$\omega$ trivial on $U_K$ and $-1$ on any uniformizer of $K$. Note
that $a(\tbts)\neq 0$. Otherwise, since $\pi_K=\pi_F \in F$ in this
case, $\tbts(\pi_K)=1$ so that $\tbts$ is trivial. Then
$\theta=\overline{\theta}$
contradicting the regularity of $\theta$. So $a(\tbts)\geq 1$.\\
We divide our counting into
mainly 3 cases:\\
{\bf Case 1:} $a(\tbts)<a(\chi)$ \\
We have $\epsilon(\chi^{-1},\psi_0)=\ot(-x_0/2)
\ot^{-1}(y_{\chi^{-1}})$ by equation (\ref{epsdeligne}). Since $\ot$
is trivial on units in the unramified case, let
$y_{\chi^{-1}}=\pi_F^{-a(\chi)}$. So we have
$\epsilon(\chi^{-1},\psi_0)=(-1)^{a(\chi)+t}$ where $t=v_F(2)$.
Since $a(\tbts)<a(\chi)$ we have
$\epsilon(\chi^{-1}\tsbt,\psi_0)=(-1)^{a(\chi)+t}=\epsilon(\chi^{-1},\psi_0)$.
So all the $\chi$'s are such that all $\chi\theta$ will occur in
${\rtheta}_+$ or all will occur in ${\rtheta}_-$ depending on whether $a(\chi)$  is even or odd.\\
{\bf Case 2:} $a(\chi)<a(\tbts)$\\
 In this
case,$a(\chi^{-1}\tsbt)=a(\tbts)$. So
$\epsilon(\chi^{-1}\tsbt,\psi_0)=(-1)^{a(\tbts)+t}$. Also,
$\epsilon(\chi^{-1}\tsbt,\psi_0)=(-1)^{a(\chi)+t}$. $\chi\theta$
will occur in $\rtheta$ if and only if $a(\chi)=a(\tbts)(mod\,2)$.\\
{\bf Case 3:} $a(\chi)=a(\tbts)$\\
Here we have two possibilities:
\begin{enumerate}
\item{$a(\chi^{-1}\tsbt)<a(\chi)$ or $a(\chi)<a(\chi^{-1}\tbts)$: In this case, if $a(\chi)=a(\tbts)(mod\,2)$  then $\chi\theta$ will occur in $\rtheta$.}
\item{$a(\chi^{-1}\tsbt)=a(\chi)$: In this case $\chi\theta$ will occur in $\rtheta$.}
\end{enumerate}
\rem Since by theorem(\ref{prasad}), $\lambda\theta$ appears in
${\rtheta}_+$ (respectively ${\rtheta}_-$) if and only if
$\lambda\theta$ does not appear in ${\rtheta'}_+$ (respectively
${\rtheta'}_-$), all the theorems proved in this paper have their
obvious $D_F^*$ analogues.\\

{\bf Acknowledgements:} I would like to thank R. Tandon, University of Hyderabad, India for some tedious discussions and some helpful suggestions and D. Prasad, Tata Institute of Fundamental Research, India for motivating me towards the problems we discussed in this exposition.


\begin{thebibliography}{abcd}


\bibitem[D]{deligne} P. Deligne, Les Constantes locales de
l'\`{e}quation fonctionelle de la fonction Ld'Artin d'une
representation orthogonale, Inv. Math. 35 (1976) 299-316.

\bibitem[K-T]{tandon}P.A. Kameswari; R. Tandon, A converse
theorem for epsilon factors, J. Number Theory 89 (2001) 308-323.

\bibitem[N-T]{vishnu} K. Vishnu Namboothiri; Rajat Tandon, Completing an extension of
Tunnell's theorem, J. Number Theory 128 (2008) 1622-1636.


\bibitem[P1]{prasad}D. Prasad, On an extension of a theorem of Tunnell,
Compositio Math. 94 (1994) 19-28.

\bibitem[P2]{prasad2}D. Prasad, Relating invariant linear
forms and local epsilon factors via global methods, with an appendix
by H. Saito, Duke J. of Math. 138(2)(2007), 233-261.

\bibitem[Ta]{tate}J. Tate, Number Theoretic Background, in Automorphic Forms,
Representations, and $L$-function, AMS Proc. Symp. Pure Math.
33(2) (1979) 3-26.
\bibitem[Tu]{tunnell} J. Tunnell, Local epsilon factors and characters of
$GL(2)$, American Journal of Math. 105 (1983) 1277-1307.

\end{thebibliography}
\end{document}